\documentclass[a4paper]{article}
\usepackage{newclude}
\usepackage{graphicx} 
\usepackage{amsmath, amsthm, amsfonts, amssymb, amscd}
\usepackage{setspace}
\usepackage{enumitem}
\usepackage{latexsym}
\usepackage{verbatim}
\usepackage[colorlinks=true,hypertexnames=false,linkcolor=blue,citecolor=blue, backref=page]{hyperref}
\usepackage{mathrsfs}
\usepackage{color} 
\usepackage[usenames,dvipsnames]{xcolor}
\usepackage{cite}
\usepackage[noblocks]{authblk}
\usepackage[utf8]{inputenc}
\usepackage{pgf,tikz}
\usepackage[pagewise]{lineno} 
\usetikzlibrary{arrows}
\pdfoutput=1

\advance \textheight by \topskip
\numberwithin{equation}{section}
\usepackage{fancyhdr}
\newtheorem{teo}{Theorem}

\newtheorem{prop}{Proposition}[section]
\newtheorem{lem}[prop]{Lemma}
\newtheorem{obser}[prop]{Remark}
\newtheorem{cor}[prop]{Corolary}
\newtheorem{defi}[prop]{Definition}

\makeatletter
\renewenvironment{proof}[1][\proofname]{\par
\pushQED{\qed}%
\normalfont \topsep6\p@\@plus6\p@\relax
\trivlist
\item\relax
{\bfseries
#1\@addpunct{.}}\hspace\labelsep\ignorespaces
}{%
\popQED\endtrivlist\@endpefalse
}
\makeatother
\def\ackname{Acknowledgements}%
\newenvironment{ack}[1][\ackname]%
{\ifx#1\empty\else\subsection*{#1.}\fi\par}
{\par}
\usepackage{cancel}
\usepackage[normalem]{ulem}
\newcommand\redst{\bgroup\markoverwith{\textcolor{red}{\rule[0.5ex]{2pt}{0.7pt}}}\ULon}
\newcommand\bluest{\bgroup\markoverwith{\textcolor{blue}{\rule[0.5ex]{2pt}{0.7pt}}}\ULon}

\usepackage[textwidth=2.7cm, color=red!20]{todonotes}
\newcounter{nota}[section]

\newcommand{\ordem}{ o }

\newcommand\restric[2]{ {
	\left.\kern-\nulldelimiterspace 
	#1 
	\vphantom{\big|} 
	\right|_{#2}
	}}

\usepackage{xparse}
\NewDocumentCommand{\hess}{om}{%
\IfNoValueTF{#1}
{{\rm Hess} (#2)} {{\rm Hess} (#2) \Big\lvert _{#1}} 
}
\usepackage{yhmath}
\usepackage{mathtools}
\newcommand{\toi}{\, \rotatebox[origin=c]{-90}{\scalebox{0.9}[1.4]{$\mathclap{{\curvearrowleft}}$}}}

\title{
On the role of the surface geometry in convex billiards}
\author{
\linespread{0.8} 
M.J. Dias Carneiro, S. Oliffson Kamphorst, S. Pinto-de-Carvalho \vskip -0.45cm
{\small Departamento de Matem\'atica, UFMG, Brasil}

C. H. Vieira Morais \vskip -0.45cm
{\small Departamento de Matem\'atica, UFES, Brasil}
}

\begin{document}
\maketitle
\begin{abstract}
\noindent 
This work presents a framework for billiards 
in convex domains on two dimensional Riemannian manifolds. 
These domains are contained in connected, simply connected open subsets which are totally normal.
In this context,  some properties that have long been known for billiards on the plane are established.
We prove the twist property of the billiard maps and establish some conditions for the existence and non-existence of rotational invariant curves. 
Although we prove that Lazutkin's and Hubacher's theorems are valid for general surfaces, 
we also find that Mather's theorem does not apply to surfaces of non-negative curvature.
\end{abstract}

\section{Introduction}

The billiard dynamics consists of the motion of a particle inside a bounded region. The particle moves freely inside the region and undergoes elastic collisions with the boundary.
Plane convex billiards where  introduced by Birkhoff \cite{birk27} in the early 1900's and since then have been  widely studied as, although of simple formulation, they exhibit very interesting dynamics ranging from integrability to ergodicity depending on the shape of the boundary.

From the flat plane, the study of billiards evolved into that of
surfaces of constant curvature \cite{veselov}. 
In this work we go further and present a framework to study convex billiards on general surfaces
{investigating  some properties of these dynamical systems which are known to hold on the flat plane  
\cite{tabach-geom,cher-mark} 
and on constant curvature surfaces \cite{luciano,florentin}.} 
{Apart from some restrictions on the surfaces and after a characterization of the concept of convexity,}
we establish that the dynamics in convex billiards on surfaces is also defined by conservative differentiable twist maps of the cylinder {which allows the use of classic tools}.
{Our study then focuses on the existence and non-existence of invariant rotational curves.  
This is an important point to characterize the global dynamics, as the existence of invariant curves prevents ergodicity while the non-existence prevents integrability.}

From a more general perspective,  
we aim to understand to what extent the {\em shape} of the surface influences the billiard dynamics.
Our goal is to determine which dynamical properties depend on the geometric properties of the boundary, such as convexity and geodesic curvature, and which depend on the geometry of the ambient surface. 

It is worthwhile to mention that the study of billiards on general surfaces is connected to the study of billiards in the flat plane with a potential. Examples of billiard motion under a gravitational \cite{lehtihet-grav,jaud-grav} or magnetic \cite{robnik,gutkin-mag,bialy-magnetic} field  have been investigated,
as well as in media with a variable index of refraction \cite{barutello-refrac,glen-refrac}.
So, billiards naturally appear in different areas of science. 
Applications such as cosmological billiards  \cite{damour-cosmo} and the asteroseismology  on rapidly rotating stars \cite{astro-spc,lignieres-star}
are also more general situations of a billiard dynamics. 
We believe that our approach can help to shed some light on these problems.

{In what follows we describe our main results and the structure of this work.} 

We begin by establishing in Section~\ref{sec:preliminares}, a class of surfaces and curves where the billiard dynamics is well defined.    
More precisely, our surfaces will be totally normal neighborhoods of two dimensional Riemannian manifolds.
In this context, we discuss the notion of convexity of curves in connection with the geodesic curvature
(Proposition~\ref{c2:prop:curv.pos.entaoconv}).
The billiard dynamics takes place inside convex regions limited by an at least $\mathscr C^2$ regular closed curve on such surfaces.  
As usual, the billiard map is defined by the successive points of impacts with the boundary and the direction of motion so that  it is given by a two dimensional map of the cylinder.

The basic results about billiards on these general surfaces are obtained in Section~\ref{sec:twist}. 
Here we  prove that the billiard map is a conservative twist map (Theorem~\ref{teo:eu.twist}).
{In particular, this allows one to use variational methods and Aubry-Mather's theory  
 to show the existence of special compact invariant subsets with {interesting dynamics}.
It follows from the Twist Property that if the boundary of the billiard is sufficiently regular, its neighborhood contains a collection of caustics, corresponding to a collection of  invariant rotational curves close to the boundary of the phase space (Lazutkin invariant curves \cite{lazutkin}).
This is the content of  Theorem~\ref{teo:eu.lazutkin}.

In the following sections we continue to address the question of the existence of invariant rotational curves.
We prove that if the boundary has points where the curvature is discontinuous, its neighborhood is caustic-free.
This result was originally proved by Hubacher \cite{hubacher} for billiards in the flat plane. 
More recently, a similar technique was used to prove that the result also holds on constant curvature surfaces \cite{florentin}.
Section~\ref{sec:hubacher} contains the proof of Hubacher's Theorem (Theorem~\ref{teo:eu.hubacher}) in our context of billiards on general surfaces.

Finally, in Section~\ref{sec:mather}, we generalize Mather's theorem \cite{mather82} (which states that if the boundary of a billiard has a point of zero curvature then there are no invariant rotational curves) to non positive curvature surfaces 
(Theorem~\ref{teo:eu.mather.K-}) and 
we establish sufficient conditions for the theorem to hold in the general case (Theorem~\ref{teo:eu.mather.K+}).
However, we also prove that these conditions are necessary, as we construct, in Section~\ref{exemplo}, an example of a constant width curve on the sphere with a point of zero curvature.
In view of  previous work on convex billiards on constant curvature surfaces, 
this  seems to be, at our knowledge,  the first result where the geometry of the surface intervenes, even though some differences were already  found in the case of polygonal billiards \cite{spina}. 

To conclude, we point out that our initial results, presented  above,  
seem to indicate that some semi local properties of the billiard dynamics
depend only on the geodesic curvature of the boundary and not on the curvature of the surface, which may however impact on global properties such as integrability. 
For instance, it is known that the (geodesic) circle is the only totally integrable convex billiard in constant curvature surfaces \cite{bialy93,bialy13}. However, there are general surfaces where the geodesic circle is not integrable  \cite{marseille}. 
The complete integrability of the circle is part of the more general  Birkhoff's conjecture, that the only integrable billiards are ellipses, which has not yet been proven even in the planar case, although several results in this direction have been recently obtained.

We would like also to observe that although  widely studied, a  the fully description of the generic properties of planar convex billiards is still missing.
By considering billiards on general surfaces, the class of perturbations to be considered is expanded and the scope of what stands for generic may be broaden. We hope to be able to identify which properties are truly persistent and also to uncover some new features.


\section{Setup: preliminaries, basic definitions and general geometric properties
}
\label{sec:preliminares}

Let $M$ be a  two dimensional ${\mathscr C}^\infty$ complete manifold with a Riemmanian metric $g$
and consider an open,  connected and simply connected subset 
$U \subset M$ such that $(U,g)$ is a totally normal neighborhood. This will be our billiard's surface.
By a {\em totally normal neighborhood} we  mean 
that  for any point $p \in U$, there is an open neighborhood $\tilde U_p$ of $0 \in T_pM$ such that $\restric{\exp_p}{\tilde U_p} : \tilde U_p \to U$ is a diffeomorphism. 
From now on, $\exp_p$ will denote this restriction
{and we will also drop the subscript $p$ when the context is clear}. 
In particular, $U$ is covered by a single chart by the exponential map and is thus naturally oriented.
The unitary disc with the hyperbolic metric is an example
of such subsets. 
So is a domain of a chart of a totally geodesic neighborhood of a point in general  Riemannian surfaces.

By restricting the class of surfaces we are going to work with, we avoid some problems. 
For instance, the very definition of the billiard motion clearly relies on the existence and uniqueness of geodesics.
As an example of the kind of ``strange'' situation which may occur in a more general surface,  let us consider
an hemisphere $\Omega_1$ of the sphere $\mathbb S^2 \subset \mathbb R^3$ and let $\gamma_1 = \partial \Omega_1$. Given $p \in \gamma_1$, any geodesic through $p$ is a circle which also intersects $\gamma_1$ at $-p$, the antipodal point of $p$. This implies that all orbits by the billiard map $T_1$ in $\Omega_1$ are 2-periodic. 
In particular,  $T_1$ isn't a twist map. 
Nevertheless, a wide variety of quite general surfaces have the properties we require.


Given two points  $p, q$ on a billiard surface  $U \subset M$, the distance between them is given by 
$d(p,q) =  |{\exp}^{-1}_p(q)| = |{\exp}^{-1}_q(p)|$, where $| \cdot |$ is the norm induced by $g$,
 i.e. $\langle u,v\rangle= g(u,v)$.
Now, let $\eta$ be the geodesic of $M$ through them which, although exists and is unique, may not be totally contained in $U$. 
By a {\em geodesic} of $U$ we mean the restrictions of $\eta$ to maximal segments in $U$. 
With this definition, geodesics do not have self intersection.    However,  the geodesic through two points may not exist.
From now on, the term {\em geodesic} will always be used in the context of $U$. 
We also assume that geodesics are normalized in the sense that $|\eta'| = 1$.
Any geodesic
separates $U$ into two distinct connected components, each one simply connected.


We will consider simple {$\mathscr C^1$} regular parametrized curves contained in $U$, meaning that the trace is contained in $U$. 
Closed curves will be positively oriented according to $\partial U$.
Since a simple regular closed curve is a Jordan regular curve contained in $U$, its complement is the union of two open sets. We then define the {\bf interior} of a simple closed curve $\gamma$ as the simply connected and bounded
component of $U \backslash \gamma$. The second component is called the {\bf exterior}.

\begin{defi}
A subset $\Omega$ such that $\bar{\Omega} \subset U$ has the {\bf convexity property} if for any two points in $\Omega$ there is a geodesic segment in $\Omega$ connecting them.
\end{defi}
If a subset has the convexity property, two points are connected by exactly one geodesic segment in it.
Moreover, it
is simply connected, it is a star set from any interior point and its closure also has the convexity property.  
\begin{defi}
An open set $\Omega \subset U$ is {\bf convex} if it has the convexity property and its boundary $\partial \Omega$ is a ${\mathscr C}^1$ regular closed curve which does not contain any geodesic segment.
The length of the geodesic segment between two points $p,q \in \bar{\Omega}$ is given by 
$d(p,q)$. 
The boundary of a convex set is a {\bf convex curve}. 
\end{defi}

It is worthwhile to observe that by our definition, convex means in fact what is usually called strictly convex in the literature.

In ${\mathbb R}^2$ a convex set is always contained in one of the half-planes determined by a tangent line. 
However, in a Riemannian surface we may not have such a global notion.
Nevertheless, this notion can be recovered for $U$, since a geodesic separates a totally normal neighborhood into two components, thus defining sides.

\begin{defi}
A regular simple curve $\gamma \subset U$ is {\bf supported} at a point $p$ if $\gamma\backslash \{p\}$ is entirely contained in one of the connected components determined by the tangent geodesic at $p$.
We say that $\gamma$ is {\bf locally supported} at $p$ if there is a ball $B$ of center $p$ such that the restriction 
$\gamma_p$ (the arc of the connected component of $\gamma \cap B$ containing $p$) is supported at $p$.
\label{def:suportada}
\end{defi}

Let  $\gamma$ be a $\mathscr C^1$ regular,
closed and simple curve with positive orientation, $\Omega$ its interior and $\Omega'$ its exterior. Given a point $p \in \gamma$, let $\eta$ be the tangent geodesic at $p$. If $\gamma$ is locally supported at $p$, let $\eta_p$ be the restriction of $\eta$ to the ball $B$ as in Definition \ref{def:suportada}. Then, either $\eta_p \cap \Omega = \emptyset$ or $\eta_p \cap \Omega' = \emptyset$.
In a flat space, the former happens if the curvature $\kappa(s)$ of $\gamma$ at $p$ is positive and the later if $\kappa(s) < 0$. If $\kappa(s) = 0$, $\gamma$ may or may not be locally supported in $p$. 
We will show that the same is true for non flat metrics. 

Throughout this work, we will alternately use  different parameterizations of a surface in order to prove our statements. Given an orthonormal and  positive basis for $T_pU$, $\{ v_1,v_2 \} $, we consider the following definitions:

\begin{defi}
The parametrization by {\bf Normal Coordinates} is $\Phi_N: V_N \subset \mathbb R^2 \to U$ given by 
$\Phi_N(x,y) = \exp_p( x v_1 + y v_2)$.
The coordinate vector fields are $X_N (\Phi_N(x,y)) = \dfrac{\partial \Phi_N}{\partial x} (x,y)$ and $Y_N (\Phi_N(x,y)) = \dfrac{\partial \Phi_N}{\partial y} (x,y)$. 
By Gauss's Lemma,  this parametrization,  {straight lines} through the origin in $\mathbb R^2$  are mapped   onto geodesics of $U$ which are orthogonal to the image of circles centered at the origin.
\end{defi}

\begin{defi} The parametrization by {\bf Polar Coordinates} is $\Phi_P: V_P \subset \mathbb R^+ \times \mathbb R \to U$ given by  $\Phi_P(r,\theta) = \exp_p r( \cos \theta v_1 + \sin \theta v_2)$.
	The coordinate vector fields are $X_P (\Phi_P(r,\theta)) = \dfrac{\partial \Phi_P}{\partial r} (r,\theta)$ and $Y_P (\Phi_P(r,\theta)) = \dfrac{\partial \Phi_P}{\partial \theta} (r,\theta)$.  $X_P$ is unitary and, for fixed $\theta$, $Y_P$ is an orthogonal Jacobi field along the geodesic $r \mapsto \Phi_P(r, \theta)$.
\end{defi}

\begin{defi}
{Let $\beta : (-\varepsilon, \varepsilon) \to U$ be a geodesic with $ \beta(0) = p$ and $\beta'(0) = v_2$}. The parametrization by {\bf Fermi Coordinates} is $\Phi_F: V_F \subset \mathbb R^2 \to U$ given by $\Phi_F(x,y) = \exp_{\beta(y)}( x \cdot v(y))$, where $v$ is a unitary, parallel vector field perpendicular to $\beta$ with $v(0) = v_1$. The coordinate vector fields are denoted by $X_F, Y_F$.
Horizontal lines are mapped  by $\Phi_F$ into geodesics of  $U$ which are orthogonal to the image of vertical lines i.e.   $X_F$ and  $Y_F$ are orthogonal.
\end{defi}
 
\begin{lem}
\label{lem:curvatura:v2}
Let  $c(t) =  (t,y(t)) $ be a ${\mathscr C}^2$ graph in $V \subset \mathbb R^2$ with 
$c(0) = (0,0)$ and $c'(0) = (1,0)$. Then the geodesic curvature $\kappa(0)$ of $c$ is equal to the geodesic curvature $\kappa_N(0)$ of  $c_N(t)  = \Phi_N(c(t))$ at $t = 0$. It is also equal to the geodesic curvature $\kappa_F(0)$ of $c_F(t) = \Phi_F(c(t))$.
\end{lem}

\begin{proof}  

Using the formula for curvature in flat spaces we get
$${\kappa = \dfrac{x'y'' - x''y'}{\left( (x')^2 + (y')^2 \right)^{3/2}} \implies \kappa(0) = y''(0)}$$
Now we proceed to compute $\kappa_i(0)$  where $i =F \hbox{ or }  N$ (for Fermi or Normal coordinates). 
On the one hand, $c_i'(t) = X_i + y'(t) Y_i$, then
\begin{eqnarray*}
\dfrac{Dc_i'}{dt} &=& \nabla_{c_i'(t)}{c_i'(t)} =  \nabla_{c_i'(t)} (X_i + y'(t) Y_i)\\
&=& \nabla_{c_i'(t)} X_i + y''(t) Y_i + \nabla_{c_i'(t)} Y_i \\
&=& \nabla _{X_i} X_i + y'(t) \nabla_{Y_i} X_i + y''(t) Y_i + \nabla_{X_i} Y_i + y'(t) \nabla_{Y_i} Y_i
\end{eqnarray*} 
{$\nabla_{X_i} X_i = \nabla_{Y_i} X_i = 0$ in $p$ since $x\mapsto \Phi_i(x,0)$ and $y\mapsto \Phi_i(0,y)$ are geodesics. Hence, for $t = 0$, $\dfrac{Dc_i'}{dt}(0)  = y''(0) Y_i$, as $y'(0) = 0$. }

On the other hand,  
$$\dfrac{d}{dt} \dfrac{1}{|c_i'(t)|} = - \dfrac{\langle \dfrac{Dc_i'}{dt}, c_i'(t) \rangle}{|c_i'(t)|^3} 
$$
which is equal to  0 for $t = 0$. 

Since $Y_i$ is the unitary normal vector of the curve $c_i$ at $t=0$ we conclude that
$$
\kappa_i(0) Y_i = \dfrac{Dc_i'}{dt}(0) \cdot \dfrac{1}{|c_i'(0)|} + c_i'(0) \cdot \left[ \dfrac{d}{dt} \dfrac{1}{\left|c_i'(t)\right|} \right]\Bigg|_{t=0}  = y''(0)Y_i$$
Hence $\kappa(0) = \kappa_N(0) =  \kappa_F(0)$.
\end{proof}

{Using a parametrization $\Phi_N$ by normal coordinates centered at $p$, by the above lemma, 
we have that $\left(D\Phi_N\right)_0 \cdot \kappa v =  \dfrac{D\gamma'}{ds}$, where $v$ is the normal unit vector to $\gamma_N = \Phi_N^{-1}(\gamma)$ at $0 = \Phi_N^{-1}(p)$ such that $\{ \gamma_1', v\}$ is a positive basis. 
As it is a well known fact that in the flat plane $v$ points toward $\gamma_N$'s interior of $\gamma_N$ has positive orientation and $v$ points to its exterior otherwise, we obtain the following equivalent result for general surfaces.
}

\begin{prop} 
\label{prop:curvatura.aponta}
Let $\gamma \subset U$ be a regular curve  and $\kappa$ its curvature at a point $p \in \gamma$. If $\kappa \neq 0$ then $\gamma$ is locally supported at $p$.
Moreover, if $\gamma$ is closed, simple and has positive orientation, $\dfrac{D\gamma'}{ds}$ points to its interior if $\kappa > 0$. Otherwise,  $\dfrac{D\gamma'}{ds}$ points to the exterior of $\gamma$. 
\end{prop}

We establish in the next proposition that curves with positive curvature are indeed convex. 
\begin{prop}
\label{c2:prop:curv.pos.entaoconv}
Let $\gamma$ be a $\mathscr C^2$ simple, closed curve contained in a 
totally normal neighborhood $(U,g)$. 

If $\gamma$ has strictly
positive geodesic curvature then it is (strictly) convex.
\end{prop}

\begin{proof}
Let $\Omega$ be the interior of $\gamma$, as given by the Jordan Curve Theorem.
Since $\gamma$ has strictly positive curvature, it cannot contain geodesic segments. 
We will prove that $\Omega$ is a convex set.

We begin by considering two arbitrary distinct points  $p, q \in \partial \Omega$ and let $\eta$ be the geodesic through them. 
Since $U$ is not necessarily a complete manifold, we have to verify that $\eta$ in fact exists. 
We can assume that $p = \gamma(0)$ to use normal coordinates centered at $p$. 
The parametrization $\Phi_N : V_N \subset \mathbb R^2 \to U$ is chosen so that the $x$ axis is mapped into the $\gamma$  tangent geodesic at $p$.  
{We can  show that a closed, simple curve with positive curvature has positive orientation so},
by Proposition~\ref{prop:curvatura.aponta},
 any geodesic segment tangent to $\gamma$ is locally outside $\Omega$
 since $\dfrac{D\gamma'}{ds}$ points towards $\Omega$. 
In particular, the $x$ axis is locally outside $\Omega_N = \Phi_N^{-1}( \Omega)$.
If $q \in \partial \Omega$ is close to $p$, by Lemma  \ref{lem:curvatura:v2}, there is a straight line segment $\eta_N$ in $V_N$  through $0 = \Phi_N^{-1}(p)$ and $ q_N  = \Phi_N^{-1}(q)$, contained in $ \overline \Omega_N$. 
As we move $q$  away from $p$ along $\partial \Omega$, it may happen that the segment $\eta_N$ is no longer contained in $\overline \Omega_{N}$. 
If this happens, there must be a position for the point $q_N$ such that  $\eta_N$ is tangent to $\partial \Omega_N$ at $q_N$ (Figure~\ref{fig:kpos.entao.convex}).
\begin{figure}[h]
\begin{center}
\includegraphics{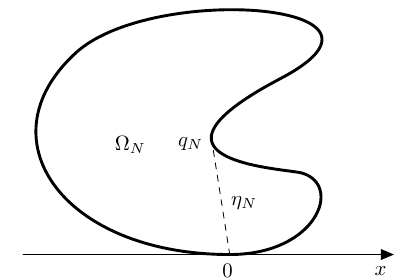}
\end{center}
\caption{Proof of Proposition \ref{c2:prop:curv.pos.entaoconv}}
\label{fig:kpos.entao.convex}
\end{figure}
In this situation, the point $q_N$  is not necessarily unique, but can be chosen such that $\eta_N \subset \overline \Omega_N$. 
However this contradicts Proposition~\ref{prop:curvatura.aponta} since $\eta  = \Phi_N(\eta_N)$ would be a geodesic segment tangent to $\partial \Omega$ at $q$ not locally outside $\Omega$.
Since $p$ was arbitrarily chosen, we have proven that there is a geodesic segment $\eta \subset \overline \Omega $ connecting any two different points in $\partial \Omega$.

To establish the convexity of $\Omega$, given two distinct points $ \tilde p, \tilde q \in \Omega$ we consider a 
parametrization $\Phi_N: V_N \to U$ by normal coordinates centered at $\tilde p$.
The straight line through $0 = \Phi_N^{-1}(\tilde p)$ and $\Phi_N^{-1}(\tilde q)$ intersects $\Phi_N^{-1}(\partial \Omega)$ at exactly two distinct points $p_N, q_N$. 
By the previous argument, the geodesic segment through $p = \Phi_N(p_N)$ and  $q = \Phi_N(q_N)$ 
is contained in $\overline \Omega$. This segment clearly contains the points $\tilde p, \tilde q$ which implies that $\Omega$ is convex. 
\end{proof}

\begin{defi}
A closed,  simple,  convex, at least $\mathscr C^2$ curve 
contained  in a totally normal neighborhood 
 will be called an {\bf oval}.
\end{defi}

The result of Proposition \ref{c2:prop:curv.pos.entaoconv} above can be extended to  piecewise $\mathscr C^2$ curves.
{Indeed, the argument in the proof relies on Proposition~\ref{prop:curvatura.aponta} which also stands using only the lateral limits of $\dfrac{D\gamma'}{ds}$ as long as they exist, are non zero and point in the same direction.}

\begin{cor}
Let $\gamma$ be a simple and closed curve contained in a totally normal neighborhood. 
Suppose that $\gamma$ is $\mathscr C^2$ and has positive geodesic curvature except at a finite number of points, where it is $\mathscr C^1$ and the lateral limits of curvature exist. 
If the curvature $\kappa$ satisfies $\kappa > \varepsilon > 0$ for some $\varepsilon$ then $\gamma$ is  (strictly) convex.
\label{cor:curv.pos.entao.conv}
\end{cor}

We conclude these preliminaries by introducing the billiard dynamics in this context.  
Let $\gamma$ be a closed,  simple,  convex, $\mathscr C^2$ curve
 i.e. an oval, contained in a totally normal neighborhood $(U,g)$. 
We will address $(U,g)$ as the {\em surface}. 
The oval  $\gamma$ is called the {\em boundary of the billiard} and 
its interior, the set $\Omega \subset U$, is called the {\em billiard table}. 
If $s$ denotes the arc length parameter of $\gamma$ and $L$ its total length,  
we consider $[0,L]$ as its domain identifying the points $0$ and $L$. 
In the sense described above $\gamma$  has a positive orientation. 

As usual, the billiard dynamics,  as  described in our introduction, is defined by the {\em billiard map}
 by associating to a departure point at the boundary and a direction of motion from it, the next impact point and the next (outgoing) direction, obtained by the reflection law. 
\begin{figure}[h]
\begin{center}
\includegraphics{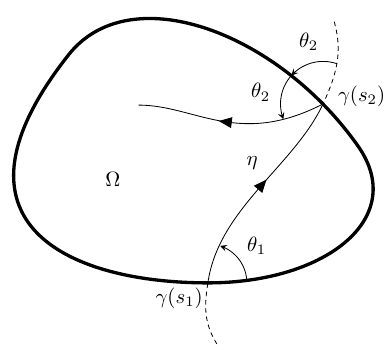}
\end{center}
\caption{Billiard variables}
\end{figure}
More precisely, given $s_1 \in [0,L]$ and $\theta_1 \in (0,\pi)$ there is a unique 
 geodesic $\eta$
such that the
oriented angle between $\gamma$ and $\eta$ at $\gamma(s_1)$ is $\theta_1$, the {\em outgoing angle}.
Since we are in a totally normal neighborhood, $\eta$ cannot be contained in a bounded set and hence must intersect $\gamma$ at points other than $s_1$. 
Therefore there is a unique $s_2 \in (0,L)$   
such that the geodesic arc through 
$\gamma(s_1),\gamma(s_2)$ is contained in $\eta$ and in 
$\Omega \, \cup \{\gamma(s_1), \gamma(s_2)\}$. 
Now,  let $\theta_2$ be the {\em incidence angle} 
measured at $\gamma(s_2)$, from 
of $\gamma$ to $\eta$.
Since the reflection angle (new outgoing angle) is equal to the incidence angle, the billard map 
$T:[0,L]\times[0,\pi] \toi$ is defined by $T(s_1,\theta_1) = (s_2, \theta_2)$ if $\theta_1 \ne 0, \pi$.  If $\theta_1 = 0$ or $\theta_1 = \pi$ we define $T(s_1, \theta_1) = (s_1,\theta_1)$.
{It is common to use the tangencial momentum $p=\cos\theta$  as variable instead of $\theta$. It is also useful to translate the billiard dynamics in terms of the evolution of Jacobi fields.}

On the other hand,
given $s_1 \ne s_2 \in [0,L)$  there is a unique geodesic $\eta$  in  $U$ from the 
$\gamma(s_1)$ to $ \gamma(s_2)$.
The oriented angle between $\gamma$ and $\eta$ at $\gamma(s_1)$ is $\Theta(s_1, s_2) = \theta_1$, the outgoing angle.


\section{Twist Property}
\label{sec:twist}

In this section we will prove our first result, the Twist Theorem (Theorem~\ref{teo:eu.twist})
and also a local version of it  (Proposition~\ref{prop:aproximacaoassintotica}).
As a consequence we prove the existence of Lazutkin curves (Theorem~\ref{teo:eu.lazutkin}) for ovals on surfaces.
We recall that, as defined in Section~\ref{sec:preliminares}, by a {\em surface} we allways mean a totally normal neighborhood of a Riemmanian surface and by an {\em oval} we mean a simple closed ${\mathscr C}^2$ curve with strictly positive geodesic curvature. 

 \begin{teo} 
{Let $\gamma$ be an oval on a surface.}
Then the billiard map on $\gamma$ is a {differentiable} twist map and preserves the measure $ds \wedge \sin \theta \, d\theta$.
\label{teo:eu.twist}
\end{teo}

The following generalization will be needed in the next section 
\begin{cor}
The result is also true if $\gamma$ is $\mathscr C^2$ except at a finite number of points  where it is $\mathscr C^1$ and the lateral limits of the curvature exist, as long as  
the curvature $\kappa$ satisfies $\kappa >  \varepsilon > 0$ for some $\varepsilon > 0$.
\label{cor:eu.twist.p.hubacher}
\end{cor}
 

To prove Theorem \ref{teo:eu.twist} we use that curves with positive curvature are indeed convex
(Proposition~\ref{c2:prop:curv.pos.entaoconv}).
For convex curves, we will show that the distance is the generating function of the billiard map.

Given a convex curve $\gamma$ in a totally normal neighborhood and $s_1 \neq s_2$, 
there is a unique geodesic segment $\eta$ through $\gamma(s_1)$ and $\gamma(s_2)$
(see \cite{manfredo-GR}, obs 3.8).
We define the distance function 
\begin{eqnarray*}
&& H: [0,L] \times [0,L] \to \mathbb R \\[2ex]
&& H(s_1, s_2) = 
\left\{ \begin{array}{lc}
d(\gamma(s_1),\gamma(s_2)), & \hbox{if $s_1 \neq s_2$;} \\
0, & \hbox{ if $s_1 = s_2$.}
\end{array} \right.
\end{eqnarray*}
where $d(\gamma(s_1), \gamma(s_2))$, as defined in the previous section, is the length of $\eta$. 
Note that $H$ is as differentiable as $\gamma$ for $s_1 \neq s_2$.  
The next  lemma shows that  $H$ is a generating function for the billiard map and concludes the proof of Theorem \ref{teo:eu.twist} together with Proposition \ref{c2:prop:curv.pos.entaoconv}.

\begin{lem} 
Let $\gamma$ be an oval on a surface
 and the function $H$ as above. For $s_1 \neq s_2$ we have
$$
\partial_1H = \dfrac{\partial H}{\partial s_1}(s_1,s_2)  = -\cos \theta_1  \qquad \textrm{ and } \qquad  \partial_2H = \dfrac{\partial H}{\partial s_2}(s_1,s_2)  = \cos \theta_2.
$$
Moreover  $\partial_{12}H = \dfrac{\partial^2H}{\partial s_1 \partial s_2}(s_1,s_2)> 0$. 
\label{lem:derivadageradora}
\end{lem}

\begin{proof}
This is a well known result for the flat plane case which has also been generalized to constant curvature surfaces \cite{??,??}.
To prove the result for a general surface,  the idea  is to notice that fixing $s_1$ and considering normal coordinates centered in $\gamma(s_1)$, we can compute the partial derivative of $H$ with respect to $s_2$ as in the flat case. 
More precisely, we fix $s_1 \in [0,L], p_1 = \gamma(s_1)$ 
and consider the parametrization $\Phi_N: V_N \subset \mathbb R^2 \to U$ centered at $p_1$, denoting
$ \gamma_N = \Phi^{-1}_{N} \, \circ \,\gamma$. 
We define the function $ h: [0,L] \to \mathbb R$ by $h (s_2) = | \gamma_N(s_2)|$.  
Since $\Phi_N$ is an isometry along lines through the origin, $H(s_1,s_2) = h(s_2)$. 
Hence
$$
\partial_2 H( s_1, s_2) = \dfrac{d}{ds_2} h(s_2) = \left \langle \dfrac{\gamma_N(s_2)}{h(s_2)}, \gamma_N'(s_2) \right \rangle
$$
By Gauss's Lemma $\partial_2H(s_1, s_2) = \left \langle v, \gamma'(s_2) \right \rangle$ where $v$ is the unitary vector tangent to the geodesic through $\gamma(s_1) $ and $\gamma(s_2)$ at $\gamma(s_2)$. 
Therefore, $\partial_2H(s_1,s_2) = \cos \theta_2$. On other hand, the symmetry of $H$  in the variables $s_1, s_2$ 
implies that 
$$
\partial_1 H(s_1, s_2) = {\partial_2 H}(s_2,s_1) = \cos (\pi - \theta_1) = - \cos \theta_1
$$
To determine $\partial_{12}H$ we write 
$$H(s_1,s_2)  \cdot \partial_1H(s_1, s_2) = - {\left \langle \gamma_N(s_2), \gamma_N'(s_1) \right \rangle}$$
Taking the derivative with respect to $s_2$ we get
$${\partial_1H \cdot \partial_2H + H \cdot \partial_{12}H =  - \left \langle \gamma_N'(s_2), \gamma_N'(s_1) \right \rangle = -  |\gamma'_N(s_2)| \cos (\theta_1 + \phi)}$$
where $\phi$ is the oriented angle from $\gamma_N(s_1)$ to $\gamma_N'(s_2)$. By our previous computation of $\partial_2 H$, we have $\cos \theta_2  = \left \langle \dfrac{\gamma_N(s_2)}{| \gamma_N(s_2)|}, \gamma_N'(s_2) \right \rangle  =  |\gamma_N'(s_2)| \cos \phi$. Therefore $ H(s_1,s_2) \partial_{12}H(s_1,s_2)  
=  |\gamma_N'(s_2)| \sin \theta_1 \sin \phi  > 0$ for $s_1 \neq s_2$.
\end{proof}

From the above lemma  we can replicate exactly the proof for the flat plane case  
 to obtain Theorem~\ref{teo:eu.twist}. 
If $\gamma$ is $\mathscr C^k $ then the billiard map is $\mathscr C^{k-1}$ in $[0,L] \times (0, \pi)$ is a result with the same proof too. We can also use the change of variables $\rho  = \cos \theta$. With this, the area form $ds \wedge d\rho$ is preserved by the billiard map.
Corollary~\ref{cor:eu.twist.p.hubacher} follows from Corollary \ref{cor:curv.pos.entao.conv} and the previous lemma. 

\begin{obser}
	As a consequence of Theorem \ref{teo:eu.twist} we can apply Birkhoff's Invariant Curve Theorem \cite{birkhoff32} for billiard maps. A {\it rotational invariant curve} by a billiard map $T: [0, L] \times [0,\pi] \toi$ is a simple, closed and not homotopically trivial curve $C\subset [0,L] \times[0,\pi]$ such that $T(C) = C$. The sets $[0,L] \times\{ 0\}$ and $[]0,L] \times \{\pi\}$ are examples of such curves. Birkhoff's Theorem states that any rotational invariant curve is a graph of a Lipschitz function.
	\label{rem:teo.birkhoff}
\end{obser}


For a convex curve in the flat plane, a linear approximation of the billiard map for small outgoing angles is obtained by replacing the curve by its curvature circle.
In the following proposition we prove that the same asymptotic expansion holds,
{although its geometric interpretation  in the context of general surfaces is not quite clear for us.}

\begin{prop}
\label{prop:aproximacaoassintotica}
Let $\gamma$ be an oval on a surface.
The associated billiard map  $T(s_1,\theta_1) = (s_2, \theta_2)$ satisfies
\begin{equation} \left\{ \begin{array}{l}
 s_2 =  s_1 + \dfrac{2 \theta_1}{\kappa(s_1)} + o\left(  \theta_1\right) \\[2ex]
 \theta_2 = \theta_1 + o\left(  \theta_1 \right) 
\end{array} \right.
\label{eq:aproximacaoassintotica}
\end{equation}
where $\ordem(x)$ denotes a function such that $\displaystyle \lim_{x \to 0} \dfrac{o(x)}{x} = 0$.
\end{prop}
\begin{proof} We can assume that $s_1 = 0$ and consider a parametrization $\Phi_N: V_N \subset \mathbb R^2 \to U$ by normal coordinates centered at $p  = \gamma(0)$. Let $\gamma_N(s) = \Phi_N \circ \gamma(s)$ and $s^N$ denote the arc length parameter for $\gamma_N$,
with $s^N(0) = 0$. 
For each $s_2$,  $s_2^N = s^N(s_2)$ is the length of $\gamma_N$ from $\gamma_N(0)$ to  $\gamma_N(s_2)$. 
The angles formed by $\gamma_N$ and the straight line segment from $0$ to $\gamma_N(s_2)$ are  
$\theta_1^N$ and $\theta_2^N$. 
Note that we have in fact a flat situation 
in $V_N$ and so the asymptotic approximation holds:
\begin{equation}
\left\{ \begin{array}{l}
 s_2^N = \dfrac{2 \theta_1^N}{\kappa_N(0)} + o\left( \theta_1^N\right) \\[3ex]
 \theta_2^N = \theta_1^N + o\left(  \theta_1^N \right) 
\end{array} \right.
\label{eq:aprox.assintotica.N}
\end{equation}
where $\kappa_N$ is the curvature of $\gamma_N$. 
Since 
the derivative of $\exp_p$ at the origin is the identity, we have
$\theta_1 = \theta_1^N$.
It follows from Lemma~\ref{lem:curvatura:v2} that 
$\dfrac{d  \theta_2}{ds}(0) = \kappa(0) = \kappa_N(0)= \dfrac{d\theta_2^N}{ds^N}(0) $
and so 
$\theta_2 = \theta_2^N + \ordem(\theta_2^N)$.
Using  the approximation of $\theta^N_2$ above (\ref{eq:aprox.assintotica.N}) we get $\theta_2 =  \theta_1 + o(\theta_1)$.

On the other hand, we notice that $\dfrac{ds^N}{ds}(0) = \left\| \dfrac{d\gamma}{ds}(0) \right\| = 1$. 
Then 
$$\dfrac{ds}{d\theta_1}(0) = \dfrac{ds}{ds^N}\cdot \dfrac{ds^N}{d\theta_1^N}  \cdot \dfrac{d\theta_1^N}{d\theta_1} = \dfrac{ds^N}{d\theta_1^N} = \dfrac{2}{\kappa_N(0)} = \dfrac{2}{\kappa(0)}$$
Hence $s_2 = \dfrac{2 \theta_1}{\kappa(0)} + \ordem (\theta_1)$.
\end{proof}

This  proposition implies in particular that the billiard map is a uniform twist for convex curves with non vanishing geodesic curvature. Moreover, the approximation \ref{eq:aproximacaoassintotica} allow us to use Douady's Corollary II-2 \cite{douady}
to prove the existence of Lazutkin's invariant curves in our setup.

\begin{teo}
{Let $\gamma$ be a $\mathscr C^7$ oval on a surface } and suppose $\gamma$ has positive geodesic curvature. Then there is a cantor set $K \subset [0, 1/2) \setminus \mathbb Q $ such that for each $\alpha \in K$ there is an rotational invariant curve $C_\alpha$ by the billiard map with rotation number $\alpha$. The restriction of the billiard map to $C_\alpha$ is conjugated to a rotation.
\label{teo:eu.lazutkin}
\end{teo}

\section{Hubacher's Theorem} \label{sec:hubacher}

In this section we extend  Hubacher's Theorem \cite{hubacher},
about the non existence of Lazutkin curves for billiards with discontinuous curvature, to general surfaces. 
This result was recently proved for constant curvature surfaces by Florentin et al. \cite{florentin} 

\begin{teo}
Let $\gamma$ be a billiard table boundary on a surface
and suppose $\gamma$ is $\mathscr C^2$ except at a finite number of points (at least one) where it is $\mathscr C^1$ and the lateral curvature  limits do exist. 
If the curvature $\kappa$ satisfies $\kappa >  \varepsilon > 0$ for some $\varepsilon > 0$ then there is a neighborhood of $[0,L] \times \{ 0 \}$ where  $[0,L] \times \{ 0 \}$ is the only rotational invariant curve for the billiard map.
\label{teo:eu.hubacher}
\end{teo}

The proof follows the original argument of  Hubacher (also used in the constant curvature case) which is based on a contradiction between
the existence of an invariant curve and Birkhoff's Theorem (remark \ref{rem:teo.birkhoff}). 
This contradiction is obtained by the construction of an orbit on the invariant curve, close to the point of discontinuity, that violates Birkhoff's ordering.
The key point is that the discontinuity of the curvature translates into a ``jump'' in the angles of a the trajectories. 

\begin{figure}[h]
	\begin{center}
		\includegraphics{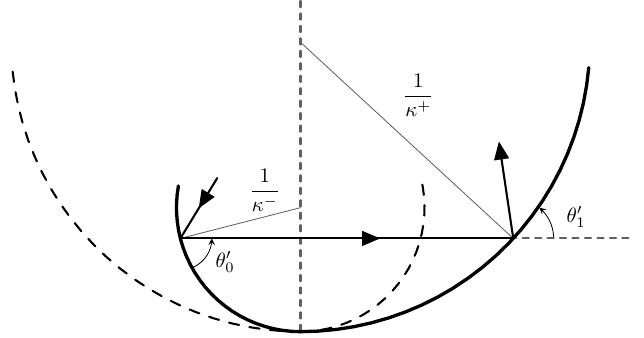}
		\caption{Two pieces of circles with different radius}
		\label{fig:hubacher-circ}
	\end{center}
\end{figure}

The original proof consists of two steps. In the first step, it is proven that any point of discontinuity of the curvature possesses a
neighborhood through which there pass no nontrivial invariant curves. 
As already mentioned, this step uses the existence of a jump in the angles of an orbit ${\cal O}(s_0',\theta_0')$ caused by the curvature discontinuity. 
This is easily verified  in the flat plane when the boundary is (locally) composed by  two pieces of circles of radius $1/{\kappa^-} < 1/{\kappa^+}$ where we have (Figure~\ref{fig:hubacher-circ})
$$ \frac{1}{\kappa^+} (1 -\cos \theta_1') =  \dfrac{1}{\kappa^-} (1 - \cos \theta_0' )
\implies \theta_1' = f(\theta_0') $$
General convex curves in the plane  are handled approximating  the curve  near the discontinuity  by its osculating circles.
The study of the derivative of $f$ provides estimates used to obtain the gap in the angles that makes the existence of invariant curves impossible 
This approach also works in constant curvature surfaces \cite{florentin} but can be tricky in a more general surface. 
We can handle this by noticing that  only information about $f'$ at the discontinuity is needed.
 This is the content of  Lemma~\ref{prop:salto} bellow, and constitutes the major difference in the proof of the theorem for non constant curvature surfaces.

\begin{lem}
\label{prop:salto}
Let $ \gamma: (-\varepsilon, \varepsilon) \to U$ be a $\mathscr C^2$ curve, except at  $s=0$ where it is $\mathscr C^1$. Assume that  the lateral  limits of the curvature at $s=0$ exist and satisfy $\kappa^- >   \kappa^+ > 0$. 
Let $\beta$ be the geodesic orthogonal to $ \gamma$ at $\gamma(0) = p$.
Given $s<0$ let $\eta$ be the geodesic through $\gamma(s)$ orthogonal to $\beta$ and $a^-, \, a^+$ the two angles
between  $\eta$ and $\gamma$. $a^-$ is orientated from $\gamma$ to $\eta$ and $a^+$, from $\eta$ to $\gamma$.
{Then $ a^+$ is a function of  $a^-$} and, for small  $ a^-$, satisfies $ a^+ < \lambda  a^-$ 
with $ 0<\lambda < 1$. 
\end{lem}

\begin{figure}[h]
\begin{center} \includegraphics{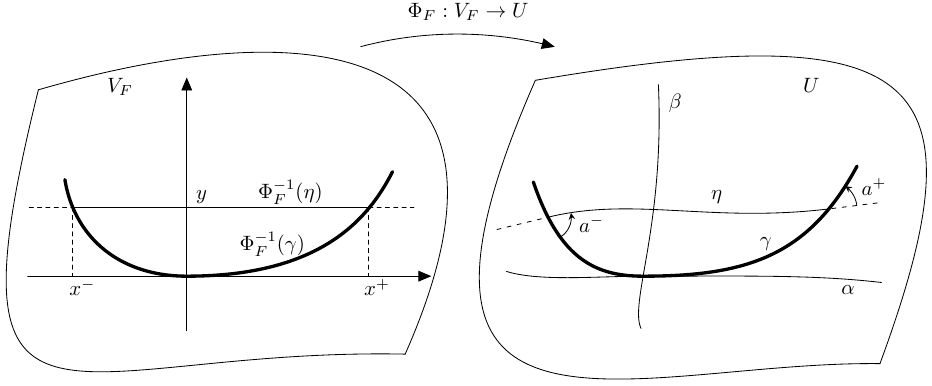}
\label{fig:hubacher}
\caption{{Lemma}~\ref{prop:salto} }
\end{center}
\end{figure}

\begin{proof}
Let $\alpha$ be the geodesic tangent to $\gamma$ at $p = \gamma(0)$ and 
$\Phi_F: V_F \subset \mathbb R^2  \to U$ the Fermi coordinates centered at $p$ such that the $x,y$ axes are mapped 
into $\alpha, \beta$ respectively. 
With this choice,  given an arbitrary $s <0$,  the geodesic $\eta$  through $\gamma(s)$ and orthogonal to $\beta$ is a coordinate curve of the parametrization (Figure~\ref{fig:hubacher}).

As $\Phi_F$ is an orthogonal parametrization the geodesic curvature of $\gamma$  at $s > 0$  is given by Liouville's  
\cite{manfredo-GD}
formula by
$$ \kappa = \dfrac{da^+}{ds} + \kappa_X  \cos a^+ + \kappa_Y \sin a^+ = \dfrac{da^+}{ds} + \kappa_Y \sin a^+ $$
where $\kappa_X,\kappa_Y$ are the  geodesic curvature of the coordinate curves at $\gamma(s)$. As $s \to 0^+$ we get 
$ \displaystyle
 \kappa^+=  \lim_{ s \to 0^+} \dfrac{da^+}{ds\ \ } $ since $\alpha$ and $\beta$ are geodesics. Analogously,  $\displaystyle  \kappa^- = \lim_{s \to 0^-} -\dfrac{da^-}{ds\ \ } $.

Since both curvatures are not zero, the functions $s \mapsto  a^-$ and $s \mapsto a^+$ are injective for $s \approx 0$. 
For each $s^- < 0$ there is a unique $x^- < 0$ such that $\Phi_F(x^-, y) = \gamma(s^-)$ for some $y > 0$. Then there is a unique $x^+ > 0$ such that $\Phi_F(x^+, y)$ is a point of $\gamma$. This defines $s^+ > 0$ by $\gamma(s^+) = \Phi_F(x^+,y)$. 
We want to estimate  the limit 
$$\lim_{a^- \to 0^-} \, \dfrac{da^+}{da^-}=
\lim_{a^- \to 0^-}  \, 
\dfrac{da^+}{ds^+} \cdot 
\dfrac{ds^+}{dx^+} \cdot 
\dfrac{dx^+}{dx^-} \cdot 
\dfrac{dx^-}{ds^-} \cdot 
\dfrac{ds^-}{da^-}$$
by analyzing each term in the chain rule product above.

Considering $r(t) = (t, y(t))$ such that $\Phi_F(r(t))$ is a reparametrization of $\gamma$, we have
$$ \lim_{x^+  \to 0^+} \, \dfrac{ds^+}{dx^+} =  \lim_{t \to 0^+} \, \dfrac{ds}{dt} =  \lim_{t \to 0^+ } \, \left| \dfrac{d}{dt} \Phi_F(r(t)) \right| = \left| \left(D\Phi_F\right)_{(0,0)} \, r'(0)\right| = 1$$
Analogously, $\lim\limits_{s^- \to 0^-} \, \dfrac{dx^-}{ds^-} = 1$. 

To analyse $\lim\limits_{x^- \to 0^-} \, \dfrac{dx^+}{dx^-}$, we use Lemma~\ref{lem:curvatura:v2} to write the Taylor expansion
$$ y(t) = \left\{ \begin{array}{rr}
\dfrac{\kappa^-}{2} t^2 + \ordem\left(t^2\right), & \textrm{ if } t < 0 \\[2ex]
\dfrac{\kappa^+}{2} t^2 + \ordem \left(t^2\right), & \textrm{ if } t > 0 
\end{array} \right.
$$
Let $ \varepsilon > 0$ such that $ 2\varepsilon < \kappa^- - \kappa^+$. For $|t|$ small enough:
$$ \left\{ \begin{array}{rr}
\dfrac{(\kappa^- - \varepsilon)}{2}t^2 < y(t) < \dfrac{(\kappa^- + \varepsilon)}{2} t^2,\,  & \textrm{ if } t < 0 \\[2ex]
\dfrac{(\kappa^+ - \varepsilon)}{2} t^2 < y(t) < \dfrac{(\kappa^+ + \varepsilon)}{2} t^2, \,  & \textrm{ if } t > 0 
\end{array} \right.   \implies 
\left\{ \begin{array}{rr}
\dfrac{ (\kappa^+ - \varepsilon)}{2} (x^+)^2 < \dfrac{(\kappa^- + \varepsilon)}{2} (x^-)^2 \\[2ex] 
\dfrac{ (\kappa^- - \varepsilon)}{2} (x^-)^2 < \dfrac{(\kappa^- + \varepsilon)}{2} (x^+)^2
\end{array}\right.
$$
This implies that $\left(\dfrac{x^+}{x^-}\right)^2$ goes to $ \dfrac{\kappa^-}{\kappa^+}$  as $\varepsilon \to 0$ and so 
$$\lim_{x^- \to 0^- } \, \left(\dfrac{dx^+}{dx^-}\right)^2 = \left( \dfrac{dx^+}{dx^-}\Bigg|_{ 0} \right)^2 = \left(\lim_{x^- \to 0^-} \dfrac{x^+ - 0}{x^- - 0}\right)^2 = \lim_{x^- \to 0^-} \left(\dfrac{x^+}{x^-}\right)^2 = \dfrac{\kappa^-}{\kappa^+}$$
Since $x^- < 0 $ and $x^+ > 0$, we get $\lim\limits_{x^- \to 0}\, \dfrac{dx^+}{dx^-} = - \sqrt{\dfrac{\kappa^-}{\kappa^+}}$. 

Finally
$$ \lim_{a^- \to 0^-} \, \dfrac{da^+}{da^-} = \kappa^+  \cdot \sqrt{\dfrac{\kappa^-}{\kappa^+}} \cdot \dfrac{1}{\kappa^-} = \sqrt{\dfrac{\kappa^+}{\kappa^-}} < 1 $$
Therefore, there is $0 <  \lambda  < 1$ such that $a^+ < \lambda a^-$ for small $a^-$.
\end{proof}

\begin{proof}[Proof of Theorem~\ref{teo:eu.hubacher}]
We follow the outline of the original proof, describing the steps leading to the result. 
To emphasize where surface curvature intervenes we subdivided the original first step into three steps as outlined bellow.
The main difference between the flat case \cite{hubacher} (and constant curvature \cite{florentin})  and ours lies on Step 1.1,
while Step 1.2 needs only  a small adaptation
and Steps 1.3 and 2 are the same. 


\label{pag:steps}
Let $\gamma$ be the billiard boundary, 
$T:[0,L] \times [0,\pi] \toi $ the associated billiard map and $s_1$ a point 
where the lateral limits of the geodesic curvature 
satisfy $\kappa^-(s_1) > \kappa^+(s_1) > 0$. 
Let us suppose 
{ that  there are nontrivial rotational invariant curves arbitrarily  close to the point $(s_1,0)$, and let $C$ such a curve. }

We will consider the lifting of $C$ in the covering space $\mathbb R \times [0,\pi]$. 
The lifting of $T$, also denoted by $T$, is fixed by taking $T(s_1,0) = (s_1,0)$.

{\bf Step 1.1:}
There is an orbit  ${\mathcal O}(s_0',\theta_0')$ belonging to the invariant curve $C$ such that  $(s_1',\theta_1') = T(s_0',\theta_0')$ satisfies 
$s_0' < s_1 < s_1'$ and $0< \theta_1' < \lambda \theta_0'$ for some $\lambda < 1$. 

{\bf (proof)} Lemma~\ref{prop:salto} above contains the construction of an orbit with required properties. 
Let $\beta$ be the geodesic  perpendicular  to the billiard boundary  $\gamma$ at the discontinuity $\gamma(s_1)$. 
By the continuity of the invariant rotational curve $C$ that we assume exists, there must be two points  
$(s_0', \theta_0')$ and $(s_1',\theta_1') = T(s_0', \theta_0')  \in C$, with $s_0' < s_1 < s_1'$, such that the geodesic segment
$\eta$ from $\gamma(s_0')$ to $\gamma(s_1')$ is perpendicular to $\beta$.
The lemma  implies that there is a neighborhood $V$ of $(s_1,0)$ and $0 < \lambda < 1$ such that 
$ 0 <  \theta_1' < \lambda \theta_0'$ if $C \cap V \neq \varnothing$.

{\bf Step 1.2:} 
Let ${\mathcal O}(s_1,\theta_1) $ be the orbit on $C$ that goes through the discontinuity and ${\mathcal O}(s_0',\theta_0')$ the orbit of Step 1.1.
There is an
integer $n$ such that 
$s_n'  < s_n$, where $T^{n-1}(s_1,\theta_1) = (s_n,\theta_n)$ and $T^n(s'_0,\theta_0') =  (s_n',\theta_n')$ .

{\bf (proof)}
Now let us consider a trajectory through the discontinuity such that the point  $(s_1, \theta_1)$ is on the invariant curve $C$.
We can choose $\delta > 0$ small enough so that $\dfrac{1 - \delta}{1+\delta} > \lambda$, which implies that
$\theta_1' <  (1+ \delta) \theta_1' <  (1-\delta) \theta_0' < \theta_0'$.
As there is a gap between $\theta_1'$ and $\theta_0'$ {as seen in this inequality}, $\theta_1$ cannot be simultaneously close to both {in the following sense: if $\theta_1 \in  \left(  \theta_1'  - \delta \theta_1' , \theta_1'  + \delta \theta_1'  \right)$, then $\theta_1 \not\in  \left(  \theta_0'  - \delta \theta_0' , \theta_0'  + \delta \theta_0'  \right)$ and vice versa.}

{
Let us assume that  $\theta_1 > (1+\delta) \theta_1'$, i.e. $\theta_1$ is ``separate" from $\theta_1'$.
As $\theta_1 > \theta_1'$ and as we are considering here small angles, the sequence $(s_k)_{k>0}$ grows faster then $(s_k')_{k>0}$. This implies that, although $s_1 < s_1'$ there is  a positive $n$ such that $s_n > s_{n}'$. The proof of this fact follow from Proposition~\ref{prop:aproximacaoassintotica}. We omit the calculations since they are the same as in \cite{hubacher,florentin}.
Analogously if $\theta_1 < (1-\delta) \theta_0'$, there is a negative $n$ such that $s_n' < s_n$.
 }

{\bf Step 1.3:}
As $s_1 < s_1'$, the inequality $ s_n' <  s_n$ of Step 1.2 contradicts the
ordering of the orbits implied by Birkhoff's theorem.
So there cannot be an invariant curve in a a neighborhood of $(s_1,0)$.

{\bf (proof)}
 Birkhoff's Invariant Curve Theorem holds for twist homeomorphisms and thus Theorem~\ref{teo:eu.twist}
	guarantees that it can be applied to convex billiards on surfaces (Remark \ref{rem:teo.birkhoff}). 
	As we noted in the beginning of the section, this theorem implies that all orbits in $C$ are ordered. On other hand, in the previous steps we have shown that there is a neighborhood $V$ of $(s_1,0)$ such that $C \cap V \neq 0$ implies the existence of non ordered orbits $\mathcal O(s_1,\theta_1), \mathcal O(s_1', \theta_1')$ of  $C$. This is a contradiction, so there can't be a nontrivial rotational invariant curve in a neighborhood of $(s_1,0)$, which concludes Step~1.

{\bf Step 2:}
Using the invariant measure, the neighborhood obtained in Step 1.3  can be extended to a whole neighborhood of $[0,L] \times \{0\}$.

{\bf (proof)}
We assume by contradiction that there is a point $(s,0)$ accumulated by invariant curves. 
Taking the limit (of a subsequence) and by Step 1, there exists an invariant curve $\tilde C$ such that $(s,0) \in C$ but $(s_1,0) \not\in \tilde C$. This, together with the Twist Property, implies the existence of an open set $W$ (Figure~\ref{fig:step2.hubacher}) such that $\mu(W) <  \mu(T(W))$, contradicting the invariance of the measure $\mu$. 
\begin{figure}[h]
\begin{center}
\includegraphics{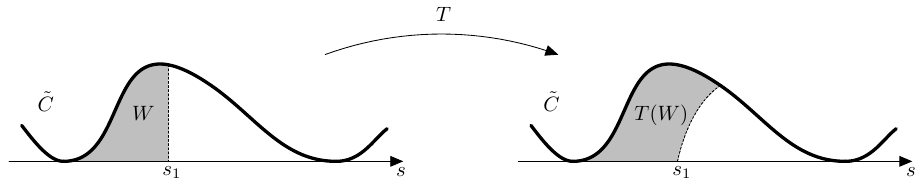}
\end{center}
\caption{Step 2 of Theorem \ref{teo:eu.hubacher}}
\label{fig:step2.hubacher}
\end{figure}
This completes the step and so the proof of the theorem.

\end{proof}

\section{Mather's Theorem} \label{sec:mather}

Mather's Theorem \cite{mather82} asserts that the existence of a point of zero curvature on a convex curve in the plane prevents the existence of rotational invariant curves for the associated biliard map.
This,  in particular, implies the existence of orbits containing points with arbitrarily small angles (i.e., glancing orbits).

We will prove that, under some hypotheses,  this results also holds for a totally normal neighborhood in a general surface. 
In fact, Mather's Theorem holds, without further assumptions, on surfaces with non-positive Gaussian curvature
(Theorem~\ref{teo:eu.mather.K-}). 
For surfaces with positive curvature, Theorem~\ref{teo:eu.mather.K+}  states that the result holds for sufficiently small curves.
Remarkably, this is not a technical hypothesis, and it is indeed possible to obtain curves with zero geodesic curvature on surfaces of positive curvature surfaces, where the billiard has an invariant rotational curve. We present an explicit example on the sphere.

\subsection{Proof of Mather's Theorem for general surfaces}

Let $\gamma$ be a $\mathscr C^2$ convex curve on
a surface (i.e. a totally normal neighborhood, as defined previously).
The associated billiard map is 
$T:[0,L] \times [0, \pi] \toi$, 
and the geodesic distance $H:[0,L] \times [0,L]  \to \mathbb R$ is the generating function as defined in Section~\ref{sec:twist}. Through this section we will also use $T$ and $H$ to denote their lifting from $[0,L]$ to $\mathbb R$.

Mather's Theorem {relies on the fact that orbits on an invariant rotational curve are minimizers of the action}.
This is expressed in the following Lemma, also due to Mather, which holds for twist maps in general, as its proof  depends only on the properties of the twist and  the generating function.

\begin{lem}[Mather]
	If $C \subset \mathbb R \times [0,\pi]$ is a rotational invariant curve for $T$,  
	then $\partial_{22}H(s_0, s_1)  + \partial_{11}H(s_1,s_2) \leq 0$ for any $(s_0,\theta_0) \in C$, where $(s_2,\theta_2) = T(s_1,\theta_1)  = T^2(s_0,\theta_0)$.
	\label{lem:inexistencia.cri}
\end{lem}

The main idea to prove the Mather's Theorem for billiards on surfaces is that we can determinate the sign of $\partial_{ii}H(s_0,s_1)$ by the analysis of a Jacobi Field along the geodesic  connecting $\gamma(s_0)$ and $ \gamma(s_1)$ as follows:
\begin{lem}
Consider the geodesic  from $\gamma(s_0)$ to $\gamma(s_1)$ and the Jacobi Field  $J:[0, H(s_0,s_1)] \to \mathbb R$ along it, satisfying $J(0)  = 0, J'(0) = 1$. 
 If the geodesic curvature of $\gamma$ is zero at $\gamma(s_1)$ then $\partial_{22}H(s_0,s_1)$ has the same sign as $J(H(s_0,s_1))J'(H(s_0,s_1))$.
	\label{lem:jacobi}
\end{lem}

\begin{proof} We fix $s_0$ and consider a parametrization $\Phi_P$ by polar coordinates around $\gamma(s_0)$. 
The map $(r, \theta) \mapsto \Phi_P(r,\theta)$ is a geodesic variation hence, for $\theta$ fixed, the vector field  
$Y_P  = \dfrac{\partial}{\partial r} \Phi_P(r,\theta)$
is a Jacobi field along the geodesic $r \mapsto \Phi_P(r, \theta)$. 
By Gauss's Lemma, $Y_P$ is an orthogonal Jacobi field. 
If $N$ is a parallel vector field along $r \mapsto \Phi_P(r, \theta)$ such that $\left \langle N, N \right \rangle  = 1$ and $\left \langle N, X_p \right \rangle  = 0$, then $Y_p = J \cdot N$ where $J$ is a function satisfying the Jacobi equation $J''(r) + K(r)J(r) = 0$ with initial conditions $J(0) = 0, J'(0) =1$ and $K$ is the Gaussian curvature along $r \mapsto \Phi_P(r,\theta)$. 
	
	Let $(r(s), \theta(s))$ be such that $\gamma(s) = \Phi_P(r(s), \theta(s))$. Our goal is to compute $r''(s_1)$ since $r(s) = H(s_0,s)$. We will omit the dependence on $s$ to make the notation lighter. Taking derivatives we get
	$$\gamma' = r' \cdot X_P + \theta' \cdot Y_P \implies \dfrac{D \gamma'}{ds} = r'' \cdot X_P + \theta'' \cdot Y_P + r' \dfrac{D}{ds} X_P + \theta' \dfrac{D}{ds} Y_p$$
	Since the geodesic curvature of $\gamma$ is zero  at $\gamma(s_1)$, the left side of the previous equation vanishes for $s = s_1$. Taking the inner product with $X_p$ and using that 
	$$\left \langle X_P , X_P \right \rangle = 1 \implies \dfrac D{ds} \left \langle X_P, X_P \right \rangle = 0 \implies \left \langle \dfrac D{ds} X_P, X_P \right \rangle = 0$$ 
we have
	$$ 0  = r'' + \theta' \left \langle \dfrac D{ds} Y_P, X_P\right \rangle $$
	On the other hand, 
	$$\dfrac{D}{ds}Y_P = r' \dfrac D{dr} Y_P + \theta' \dfrac D{d\theta} Y_P \implies r'' = - \theta' r'
	\left \langle \dfrac D{dr} Y_P, X_P \right \rangle - \theta'^2 \left \langle \dfrac D{d\theta} Y_P, X_P\right \rangle$$
	The first inner product is zero as
	$$\dfrac D{dr} Y_P = \dfrac D{dr} \dfrac{\partial}{\partial \theta} \Phi_P = \dfrac{D}{d\theta} \dfrac{\partial}{\partial r} \Phi_P  = \dfrac{D}{d\theta} X_p \implies  \left \langle \dfrac D{dr} Y_P, X_P \right \rangle = \left \langle \dfrac D{d\theta} X_P, X_P\right \rangle = 0$$
	For the second inner product we have
	$$ \left \langle \dfrac{D}{d\theta} Y_P, X_P \right \rangle  = \dfrac d{d\theta} \left \langle \dfrac{\partial}{\partial \theta}\Phi_P, \dfrac \partial {\partial r} \Phi_P \right \rangle -  \left \langle \dfrac{\partial }{\partial \theta} \Phi_P, \dfrac D{d\theta } \dfrac \partial {\partial r} \Phi_P \right \rangle  = - \left \langle Y_P, \dfrac D{dr} \dfrac \partial {\partial \theta }\Phi_P \right \rangle  = - \left \langle Y_P, \dfrac D{dr} Y_P \right \rangle $$
As $Y_P = J \cdot N$ 
we have
	$$ \left \langle Y_P, \dfrac D{dr} Y_P \right \rangle = \left \langle J \cdot N, J' \cdot N  \right \rangle  = J J'$$
	combining this with the previous equations we get $ r''(s_1) = \theta'(s_1)^2 J(r(s_1)) J'(r(s_1))$ which can be rewritten  as 
	$$ \partial_{22}H(s_0,s_1) = \theta'(s_1)^2 J(H(s_0, s_1))J'(H(s_0,s_1)) $$
	concluding the proof.
\end{proof}

{The key observation is that as a Jacobi field cannot vanish twice in a totally normal neighborhood, 
we must have $J(H(s_0,s_1)) > 0$  and the sign of $\partial_{ii}H(s_0,s_1)$ is completely determined by the sign of $J'(H(s_0,s_1))$.  In non positive curvature surfaces, we can easily control  $J'$ using Jacobi's equation.}

\begin{teo}  Let $\gamma$ be a $\mathscr C^2$ convex curve on a surface (totally normal neighborhood) of non positive Gaussian curvature. If the geodesic curvature of $\gamma$ vanishes at some point then the billiard map associated to $\gamma$ has no rotational invariant curves  besides the top and bottom boundaries of the phase space.
\label{teo:eu.mather.K-}
\end{teo}

\begin{proof}
Let $s_1$ the point where the curvature vanishes and let $s_0, s_2$ be such that 
$T(s_0,\theta) = (s_1, \theta_1), T(s_1, \theta_1) = (s_2, \theta_2)$ 
for some $\theta_0, \theta_1, \theta_2\neq 0$. 

We consider the Jacobi field $(J(0)=0, J'(0)=1)$ along the geodesic from $\gamma(s_0)$ to $\gamma(s_1)$ as in 
Lemma~\ref{lem:jacobi}. 
As $\gamma$ is contained in a totally normal neighborhood, we have $J(r) >0$ for $r \in (0, H(s_0, s_1)]$ and 
since  $J''(r) = - K(r)J(r)$, $K(r) \le 0$ implies that $J''(r) \ge 0$ and so $J'(r) \ge 1$. 
In particular $J(H(s_0,s_1)) > 0$ and $J'(H(s_0,s_1)) >0$ and it follows from Lemma~\ref{lem:jacobi} that  
$\partial_{22}H(s_0,s_1) > 0$. 

Since $\partial_{11}H(s_1,s_2) = \partial_{22}H(s_2,s_1)$,
 a similar argument shows that $\partial_{11}H(s_1,s_2)>0$ and, as in the flat plane situation,
 the theorem follows directly from 
Mather's Lemma~\ref{lem:inexistencia.cri}.
\end{proof}

\begin{obser}
By definition \cite{sullivan} Theorem~\ref{teo:eu.mather.K-} holds in the more general case of surfaces without focal points.
\end{obser}
{If the Gaussian Curvature of the surface is positive, we can no longer guarantee that $J'$ remains positive, since $J''$ may be negative for $r>0$.
However,  by continuity, as $J'(0) =1$ we will have $J'>0$  for $r$ sufficiently small. 
So, in a billiard where the distances are "small", Mather's Theorem will hold.}

\begin{teo}\label{teo:eu.mather.K+}
	Let $U$ be a totally normal neighborhood in a Riemmanian surface.
There is an open set $V \subset U$ and a constant $D> 0$ such that if $\gamma$ is a $\mathscr C^2$ convex curve contained in $V$ with diameter at most $D$ and $\kappa(s_1) = 0$ for some $s_1$, then its billiard map has no invariant curve besides the top and bottom lines.

\end{teo}

\begin{proof}
	Let $q \in U$ and consider the map $F:T{U} \to {U} \times {U}$ given by $F(q,v) = (q, \exp_q v)$. Since $(d\exp_q)_0 = \textrm{Id}$, the derivative of $F$ at $(q,0)$ is $\begin{bmatrix}
		\textrm{Id}&   0\\
		\textrm{Id}  & \textrm{Id} 
	\end{bmatrix}$. 
Then given $ M > 1$ 
there is $V \subset U$ and $\varepsilon >0$ such that $\left \| \left( d \exp_p \right)_v \right \| < M$ for any $p \in V$, $v \in T_p {U}, |v| < \varepsilon$. This implies that the solution $J$ to $J''(r) + K(r)J(r) = 0$ as stated in Lemma \ref{lem:jacobi} satisfies $J(r) < rM$ if the geodesic segment from $\gamma(s_0)$ to $\gamma(s_1)$ is contained in $V$. Reducing $V$ is necessary, we can assume that the Gaussian curvature $K$ is limited in $V$. 
Let $\tilde K >0$ an upper bound and $D < \dfrac 1{\sqrt {KM}}$. 
	
	The pair $V, D$ has the required properties. To see this, let $s_0, s_2$ such that $T(s_0,\theta) = (s_1, \theta_1), T(s_1, \theta_1) = (s_2, \theta_2)$ for some $\theta_0, \theta_1, \theta_2\neq 0$. As in the proof of the previous theorem, we only need to show that $J(H(s_0,s_1))J'(H(s_0,s_1))$. Using that $J(r) < rM$ and $K(r) < \tilde K$ we get $J''(r) = - K(r)J(r) > - \tilde K r M$. By the Mean Value Theorem and using $r, H(s_0,s_1) < D$:
	$$J'(H(s_0,s_1)) - 1  > - H(s_0,s_1) K r M \implies J'(H(s_0,s_1)) > 1 - D^2 KM > 0$$
	As $V$ is a totally normal neighborhood itself, $J$ is positive. Hence $\partial_{22}H(s_0,s_1) > 0$ by Lemma \ref{lem:jacobi}.
\end{proof}


\subsection{A ``counterexample" to Mather's Theorem} \label{exemplo}

{By controlling the sign of the {derivative of the} solutions of the Jacobi's equation, reducing the size of the neighborhood if necessary, we have managed to obtain a sufficient condition to Mather's Theorem on general surfaces.  
However, in a surface with positive curvature, it is in fact possible that this property of the Jacobi fields is not satisfied.
It is then natural to ask if, in this case, it is possible to construct a billiard table with a point of zero geodesic curvature, such that the map has an invariant rotational curve. 
The answer is yes, as we will exhibit in the following example. 
This evidences that the hypothesis in Theorem~\ref{teo:eu.mather.K+} not only technical but is indeed necessary. }

The example is a constant width curve with a point of zero geodesic curvature on the sphere. 
The line $\theta = \pi/2$ is invariant as every point on it has period two.

This construction we use to obtain a curve of constant width is quite standard in ${\mathbb R}^2$: Starting from a smooth convex arc $\alpha$ with parallel tangents at the end points and distance $\lambda$ between them equal to the sum of the curvature radius,
a closed curve is obtained by traversing from each point a distance $\lambda$  in the normal  direction. 
If the curvature of $\alpha$ is bounded from bellow by $1/\lambda$, the resulting curve is simple and convex and has constant width 
$\lambda$. 
On the sphere, however,  it is possible to obtain a smooth convex curve with a point of zero curvature.

\begin{lem}\label{lem:largura}
Consider $\mathbb S^2 =  \{ (x,y,z) \in \mathbb R^3| \ x^2 + y^2 + z^2 = 1\}$ 
the 2 dimensional unit sphere. $\mathbb S^2$ is a Riemann manifold with the induced inner product from $\mathbb R^3$. Any open hemisphere is a totally normal neighborhood of $\mathbb S^2$. 

Let $\alpha: [-L,L] \to {\mathbb S^2}$
 be a $\mathscr C^3$ arc parametrized by arc length and with non negative geodesic curvature 
$\kappa$, such that 
\begin{enumerate}[noitemsep,topsep=0pt]
\item
$\alpha'(-L) = - \alpha'(L)$
\item
 $d(\alpha(-L), \alpha(L)) = \lambda \in  \left(\pi/2, \pi \right)$ 
\item
 $\kappa(L) = \kappa(-L) = \cot (\frac \lambda 2)$
\item
  $\kappa(0) = 0$
\item
{$\kappa(s) \in [0,-\tan \lambda]$}
\end{enumerate}
If $N_\alpha$ denotes the normal vector of $\alpha$,  let $\beta : [-L,L] \to  {\mathbb S^2}$ be 
$$ \beta(s) = \alpha(s) \cos \lambda + N_\alpha(s) \sin \lambda$$
Then
$\gamma : [-L,3L] \to  {\mathbb S^2}$ defined by
$$ \gamma(s) = \left\{ \begin{array}{ll}
	\alpha(s)  &\textrm{ if } -L \leq s \leq L \\
	\beta(s-2L) &\textrm{ if } L < s \leq 3L 
\end{array}\right. $$
is a $\mathscr C^2$ convex constant width curve on the sphere. 
\end{lem}

\begin{proof}
If $p \in {\mathbb S^2}$ and $v \in T_p {\mathbb S^2}$, then $t \mapsto p \cos t + v \sin t$ is a geodesic for $ {\mathbb S^2}$. 
By the definition of $\beta$, we have $\beta(L) = \alpha(-L)$, $\beta(-L) = \alpha(L)$ and for $s\in [-L,L]$,  
$d(\alpha(s), \beta(s)) = \lambda$.

\begin{figure}[h]
	\centering
	\includegraphics{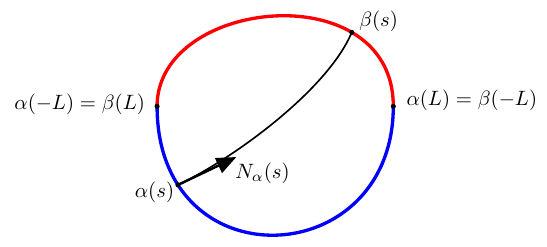}
\caption{Lemma~\ref{lem:largura}.}
\end{figure}

We have
$\beta'(s) = \alpha'(s) \left( \cos \lambda -  \kappa(s) \sin \lambda \right) $.
As  $\kappa \geq 0 $ and $ \dfrac \pi2 < \lambda <  \pi$,  $(\cos \lambda  - \kappa(s) \sin \lambda)  < 0$,
so $\beta$ is  a regular curve and  the tangent vector 
$\frac{\beta'(s)}{|\beta'(s)|} = -\alpha'(s)$.
Moreover $\beta'(s)$ 
is orthogonal  to $ (- \alpha(s)  \sin \lambda + N_\alpha(s) \cos \lambda) $.
Thus, the geodesic $t\mapsto (\alpha(s) \cos t + N_\alpha(s) \sin t)$  is orthogonal to both $\alpha$ and $\beta$ at the intersection points. It also follows that $N_{\beta}(s) =  ( \alpha(s)  \sin \lambda - N_\alpha(s) \cos \lambda) $

Computing the second derivative of  $\beta$  and using that $\langle \beta'', N_{\beta} \rangle = |\beta'|^2 \kappa_{\beta}$, we obtain
$$ \kappa_\beta(s) = - \dfrac{\sin \lambda + \kappa(s) \cos \lambda}{\cos \lambda - \kappa(s) \sin \lambda}$$
and so $\kappa_{\beta} > 0 $ for $\kappa < -\tan \lambda$.
As $\kappa(-L) = \kappa(L) = \cot(\frac \lambda 2) $ we have $\kappa_\beta(-L) =\kappa_\beta(L) = \cot( \frac \lambda 2) $ and 
$\gamma$ can be re-parametrized as a $\mathscr C^2$ regular convex curve.
\end{proof}
If $\gamma$ is a constant width curve, as in Lemma~\ref{lem:largura} above, the orthogonality property of  the geodesic from $\alpha(s)$ to $\beta(s)$ implies that the billiard trajectory between these two points has period 2 and thus the line $\theta = \pi/2$ is  an invariant rotational curve of the billiard map $\gamma$.

We conclude this section by showing that it is possible to find an arc $\alpha$ with the required properties, by exhibiting an explicit example. 

Using the standard parametrization of the sphere 
$$\Phi(\theta, \varphi) = (\cos \theta \sin \varphi, \sin \theta \sin \varphi, \cos \varphi)$$
we define an arc $ \alpha: [-1,1] \to {\mathbb S^2}$ by $ \alpha (t) = { \Phi( {\tilde \alpha(t)})}$ with 
$\tilde \alpha(t) = (\theta(t), \varphi(t))$. The functions $\theta$ and $\varphi$ are defined for $t \in [0,1]$ and 
$\mathscr C^2$ 
extended to $[-1,1]$ by $\theta(-t) = -\theta(t)$ and $\varphi(-t) = \varphi(t)$.
The idea is to translate properties 1--5 in Lemma~\ref{lem:largura} to the functions $\theta$ and $\varphi$.
To do so, we will use that
\begin{equation} \label{eqn:kappa}
\kappa v^3 + \tilde\kappa \tilde v^3 \sin \varphi = \theta' \cos \varphi \,( \theta'^2 \sin^2 \varphi + 2 \varphi'^2)
\end{equation}
where $v = |\alpha'| = \sqrt{\varphi'^2 + \theta'^2 \sin^2 \varphi}$, $\tilde v = | \tilde \alpha'| = \sqrt{\theta'^2 + \varphi^2}$ and $\kappa$, $\tilde \kappa$ are the geodesic curvature of $\alpha$ and $\tilde \alpha$.

We start by fixing the extremities by $\theta(0) = 0, \theta(1) = \pi/2$ and $\varphi(0) = \varphi_0, \varphi(1) = \varphi_1$ 
(with $\pi/2 > \varphi_0 > \varphi_1 > \pi/4$) and imposing that $\alpha$ is tangent to the parallel $\varphi=\varphi_1$ at its end points, and to the parallel $\varphi=\varphi_0$ at its middle point, i.e
\begin{equation} \label{eqn:p1}
\theta'(1) \, ,  \, \theta'(0) \ne 0, \ \ \ \varphi'(1) = \varphi'(0) = 0
\end{equation}
With this choice,  not only $\alpha'(-1) = \alpha'(1)$ 
but we also have that $\pi/2 < d(\alpha(1),\alpha(-1)) = 2 \varphi_1 = \lambda < \pi $. 

Using equation~\ref{eqn:kappa} at $t=0$ we have that 
\begin{equation} \label{eqn:p4}
\varphi''(0) = \theta'(0)  \cos \varphi_1 \sin \varphi_1 \implies \kappa (0) = 0
\end{equation} 
It is easy to verify that 
$\tilde\kappa(1) = 0$ implies that  $\kappa(1) = \cot \varphi_1$ 
so we have 
\begin{equation} \label{eqn:p3}
\varphi''(1) = 0 \implies \kappa(1) = \cot \lambda/2 
\end{equation} 

Finally, to assure that the curve is $C^{3}$ we want $\tilde\kappa'(0) = 0$ which implies
\begin{equation}\label{eqn:suave}
\varphi'''(0) = 3 \theta'(0) \theta''(0) \cos \varphi_0 \sin \varphi_0 
\end{equation}

Conditions \ref{eqn:p1}--\ref{eqn:suave} imply that the resulting curve $\alpha$ is $C^3$ and satisfies properties 1--4 in 
Lemma~\ref{lem:largura}. We notice that these conditions involve up to third derivatives of  $\theta$ and $\varphi$ at $t=0$ and $t=1$ so it is natural to use 
$$(\theta(t),\varphi(t))= \sum_{k=0}^6 (a_k,b_k) \binom 6k t^k(1-t)^{6-k} $$
 and adjust the coefficients $(a_k,b_k)$. 
This was done using a Bezier construction with GeoGebra software and the resulting curve is displayed on 
Figure~\ref{fig:exemplo-mather}.

\begin{figure}[h]
\centering
{\includegraphics{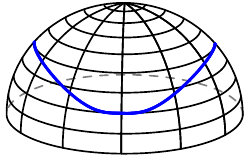}}%
\hspace*{1cm}%
{\includegraphics{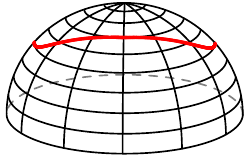}}%
\hspace*{1cm}%
{\includegraphics{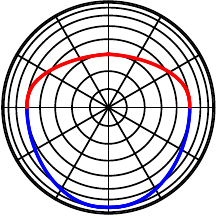}}%
\caption{A constant width curve on the sphere}\label{fig:exemplo-mather}
\end{figure}

Once the curve was obtained, Property~5 of Lemma~\ref{lem:largura} was explicitly verified. The geodesic curvature of 
the constant width curve is plotted on Figure~\ref{fig:curvatura-mather}.
\begin{figure}[h]
\centering
\includegraphics[width=0.4\linewidth]{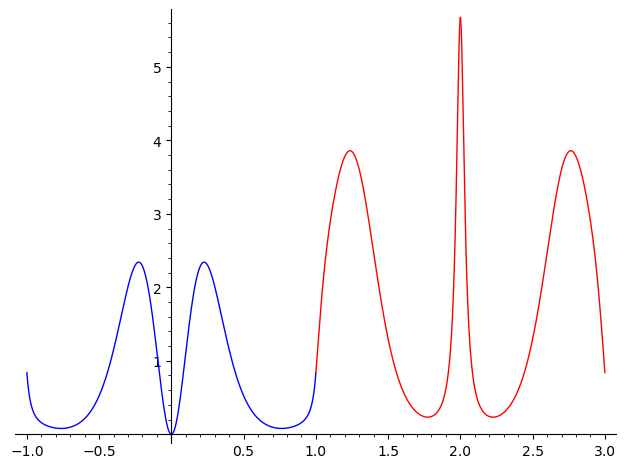}
\caption{Geodesic curvature of the constant width curve.}\label{fig:curvatura-mather}
\end{figure}

\begin{ack}
This work originated from the PhD thesis of C.H. Vieira Morais  \cite{cassiotese} which
was financed by the Coordena\c c\~ao de Aperfei\c coamento de Pessoal de N\'\i vel Superior (CAPES),  Brasil.
\end{ack}

{\linespread {1.2}
{\bf Competing interest statement declared by the corresponding author on behalf of all authors:}
The authors M.J. Dias Carneiro, S. Oliffson Kamphorst and S. Pinto-de-Carvalho have not  received payments or travel funding for the present work.  C.H. Vieira Morais was financed during his Ph.D. by a scholarship from Coordena\c c\~ao de Aperfei\c coamento de Pessoal de N\'\i vel Superior (CAPES),  Brasil from August 2016 to Mars 2021. Total funding: approximately US\$ 23,000.
}
\bibliographystyle{plain}
\bibliography{referencias}
 \end{document}